\documentclass[12pt]{article}
\usepackage{graphicx}
\begin{document}

\title{The Number of Point-Splitting Circles}
\author{Federico Ardila M.}
\date{August 3, 2001}
\maketitle

\begin{abstract}
Let $S$ be a set of $2n+1$ points in the plane such that no three
are collinear and no four are concyclic. A circle will be called
{\it point-splitting} if it has 3 points of $S$ on its
circumference, $n-1$ points in its interior and $n-1$ in its
exterior. We show the surprising property that $S$ always has
exactly $n^2$ point-splitting circles, and prove a more general
result.
\end{abstract}

\medskip

\section{Introduction}

Our starting point is the following problem, which first appeared
in the 1962 Chinese Mathematical Olympiad \cite{4}.

\medskip

{\bf Problem 1.1.} Let $S$ be a set of $2n+1$ points in the plane
such that no three are collinear and no four are concyclic. Prove
that there exists a circle which has 3 points of $S$ on its
circumference, $n-1$ points in its interior, and $n-1$ in its
exterior.

\medskip

Following \cite[p.48]{3}, we call such a circle {\it
point-splitting} for the given set of points. For the rest of
sections 1 and 2, $S$ denotes an arbitrary set of $2n+1$ points in
general position in the plane, where $n$ is a fixed integer.

There are several solutions to problem 1.1. Perhaps the easiest
one is the following. Let $A$ and $B$ be two consecutive vertices
of the convex hull of $S$. We claim that some circle going through
$A$ and $B$ is point-splitting. All circles through $A$ and $B$
have their centers on the perpendicular bisector $\ell$ of the
segment $AB$. Pick a point $O$ on $\ell$ which lies on the same
side of $AB$ as $S$, and is so far away from $AB$ that the circle
$\Gamma$ with center $O$ and going through $A$ and $B$ completely
contains $S$. This can clearly be done. Now slowly ``push" $O$
along $\ell$, moving it towards $AB$. The circle $\Gamma$ will
change continuously with $O$. As we do this, $\Gamma$ will stop
containing some points of $S$. In fact, it will lose the points
of $S$ one at a time: if it lost $P$ and $Q$ simultaneously, then
points $P,Q,A$ and $B$ would be concyclic. We can move $O$ so far
away past $AB$ that, in the end, the circle will not contain any
points of $S$.

Originally, $\Gamma$ contained all the points of $S$. Now, as it
loses one point of $S$ at a time, we can decide how many points
we want it to contain. In particular, if we stop moving $O$ when
the circle is about to lose the $n$-th point $P$ of $S$, then the
resulting $\Gamma$ will be point-splitting: it will have $A$, $B$
and $P$ on its circumference, $n-1$ points inside it, and $n-1$
outside it, as shown in Figure 1.

\begin{figure}
\centering
\includegraphics[width=2.5in]{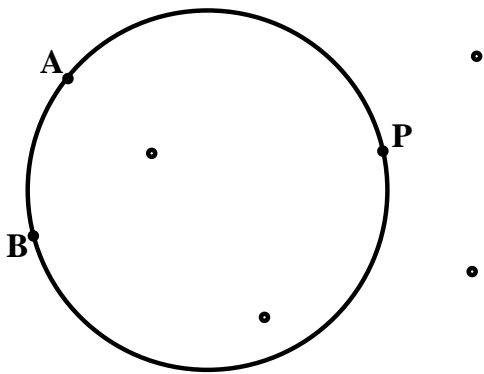}
\caption{A point-splitting circle through $A, B$ and $P$}
\end{figure}

The above proof hints that any set $S$ has several different
point-splitting circles. We can certainly construct one for each
pair of consecutive vertices of the convex hull of $S$. In fact,
the argument above can be modified to show that, for {\bf any}
two points of $S$ we can find a point-splitting circle going
through them. The reader might find it instructive to work out a
proof.

This suggests that we ask the following question. What can we say
about the number $N_S$ of point-splitting circles of $S$? At
first sight, it seems that we really cannot say very much about
this number. Point-splitting circles seem hard to ``control", and
harder to count.

We should be able to find upper and lower bounds for $N_S$ in
terms of $n$. Right away we know that $N_S \geq n(2n+1)/3$, since
we can find a point-splitting circle for each pair of points of
$S$, and each such circle is counted by three different pairs.
Computing an upper bound seems more difficult. If we fix points
$A$ and $B$ of $S$, it is indeed possible that all $2n-1$ circles
through $A$, $B$ and another point of $S$ are point-splitting.
The reader is invited to check this. This is not likely to happen
very often in a set $S$, and we can get {\it some} upper bound out
of this. However, it is hard to make this precise and get a
non-trivial upper bound.

When $S$ consists of 5 points, the situation is simple enough that
we can actually show that $N_S=4$ always. This was done in
\cite{2}. It was also proposed, but not chosen, as a problem for
the 1999 International Mathematical Olympiad. Notice that our
lower bound above gives $N_S \geq 4$.

In a different direction, problem 5 of the 1998 Asian-Pacific
Mathematical Olympiad, proposed by the author, stated the
following.

\medskip

{\bf Proposition 1.2.} $N_S$ has the same parity as $n$.

\medskip

This result follows easily from the nontrivial observation that,
for any $A$ and $B$ in $S$, the number of point-splitting circles
that go through $A$ and $B$ is odd.

The following result brings together the above considerations.

\medskip

{\bf Theorem 1.3.} Any set $S$ of $2n+1$ points in the plane in
general position has exactly $n^2$ point-splitting circles.

\medskip

Theorem 1.3 is the main result of this paper. In section 2 we
prove that every set of $2n+1$ points in the plane in general
position has the same number of point-splitting circles. In
section 3 we prove that this number is exactly $n^2$. In section
4 we present some questions that arise from our work.

\medskip

\section{$N_S$ is Constant}

\medskip

At this point, we could go ahead and prove the very
counterintuitive Theorem 1.3, suppressing the motivation behind
its discovery. With the risk of making the argument seem longer,
we believe that it is worthwhile to present a natural way of
realizing and proving that the number of point-splitting circles
of $S$ depends only on $n$. Therefore, we ask the reader to
forget momentarily the punchline of this article.

Suppose that we are trying to find out whatever we can about the
number $N_S$. As mentioned in Section 1, this number does not
seem very tractable and it is not clear how much we can say about
it. Being optimistic, we can hope to be able to answer the
following two questions.

\medskip

{\bf Question 2.1}. What are the sharp lower and upper bounds $m
= m_{2n+1}$ and $M = M_{2n+1}$ for $N_S$?

\medskip

{\bf Question 2.2}. What are all the values that $N_S$ takes in
the interval $[m, M]$?

\medskip

Question 2.1 seems considerably difficult. To answer it
completely, we would first need to prove an inequality $m \leq N_S
\leq M$, and then construct suitable sets $S_{min}$ and $S_{max}$
which achieve these bounds. To see how difficult this is, the
reader is invited to try to construct {\it any} set $S$ of $2n+1$
points for which the number $N_S$ can be easily computed.

At this point question 2.1 seems very hard, so let us focus on
Question 2.2 instead. Here is a first approach.

Intuitively, since the set $S$ can be transformed continuously,
we should expect the value of $N_S$ to change ``continuously" with
it. Suppose we start with the set $S_{min}$ (with $N_S=m$) and
move its points continuously so that we end up with $S_{max}$
(with $N_S=M$). The value of $N_S$ should change ``continuously"
as $S$ changes continuously. By ``continuity" we would guess that
$N_S$ sweeps all the integers between $m$ and $M$ as $S$ changes
from $S_{min}$ to $S_{max}$.

Right away, we know that this is not entirely true. By
Proposition 1.2 we know that the parity of $N_S$ is determined by
$n$, so $N_S$ will not sweep {\bf all} the integers between $m$
and $M$. This is not too surprising, since we haven't made
precise the meaning of the statement that the value of $N_S$
should change ``continuously" as $S$ changes continuously. The
above guess assumed that the value of $N_S$ can only jump by 1 as
$S$ is transformed continuously. (That is, if we have a set with
$k-1$ point-splitting circles and we deform it continuously into
a set with $k+1$ point-splitting circles, then somewhere in the
middle we must have had a set with $k$ point-splitting circles.)
We have no reason to assume that.

We can still hope that, as $S$ changes, $N_S$ sweeps all the
integers {\bf of the right parity} between $m$ and $M$. To show
this, we would have to show that the value of $N_S$ can only jump
by 2 as $S$ is transformed continuously. This is a reasonable
statement which we can try to prove.

In any case, the natural question to ask is what kind of
``continuity" the value of $N_S$ satisfies as $S$ changes
continuously. We certainly expect that if two sets $S$ and $T$
look very very much alike, then the difference $N_S-N_T$ should
be small. We have to find a way to make this statement precise.

\medskip

Suppose we have sets $S_{min}=\{P_1, \ldots, P_{2n+1}\}$ and
$S_{max}=\{Q_1, \ldots, Q_{2n+1}\}$ that achieve the upper and
lower bounds for $N_S$, respectively. Now slowly transform
$S_{min}$ into $S_{max}$: first send $P_1$ to $Q_1$ continously
along some path, then send $P_2$ to $Q_2$ continuously along some
other path, and so on. We can think of our set $S$ as changing
with time. At the initial time $t=0$, our set is $S(0)=S_{min}$.
At the final time $t=T$, our set is $S(T)=S_{max}$. In between,
$S(t)$ varies continously with respect to $t$. How does
$N_{S(t)}$ vary ``continuously" with time? How small can we make
$N_{S(t + \Delta t)} - N_{S(t)}$ for small enough $\Delta t$?
This is the question we need to ask.

\medskip

{\bf Technical Remark.} As we move from $S(0)$ to $S(T)$
continuously, it is likely that several intermediate sets $S(t)$,
with $0 < t < T$, are not in general position. Strictly speaking,
we should only consider those times $t$ when $S(t)$ is in general
position; when $S(t)$ is not in general position, we should
decree that $S(t)$ is undefined, and have a discontinuity at $t$.

We shall see that we can go from $S(0)$ to $S(T)$ with only
finitely many such discontinuous points. At such a discontinuity
$t$, we still need to know how small we can make $N_{S(t + \Delta
t)} - N_{S(t - \Delta t)}$ for small $\Delta t$.

\medskip

For small enough $\Delta t$, the set $S(t + \Delta t)$ is a very
slight deformation of $S(t)$. What is missing is an understanding
of what can make $N_S$ change as the set $S$ changes very
slightly from $S(t)$ to $S(t+\Delta t)$, and how small this
change is. Let us answer this question.

Notice that, in the way we defined the deformation from $S_{min}$
to $S_{max}$, the points of $S$ moved only one at a time. Let us
focus for now on the interval of time where $P_1$ moves towards
$Q_1$.

Suppose that the number $N_S$ changes between time $t$ and time
$t + \Delta t$. Then it must be the case that for some $i,j,k$
and $l$ the circle $P_iP_jP_k$ contained (or did not contain)
point $P_l$ at time $t$, but at time $t + \Delta t$ it does not
(or does) contain it. For this to be true, it must have happened
that somewhere between times $t$ and $t + \Delta t$, these four
points must have been concyclic, or three of them must have been
collinear. Since $P_1$ is the only point that has moved, we can
conclude that $P_1$ must have crossed a circle or a line
determined by the other points; this is what caused $N_S$ to
change. We will call the circles and lines determined by the
points $P_2, P_3, \ldots, P_{2n+1}$ the {\it boundaries}.

We can choose the path along which $P_1$ is going to move towards
$Q_1$. To make things easier, we may assume that $P_1$ never
crosses two of the boundaries at the same time. This can clearly
be guaranteed: we know that these boundaries intersect pairwise
in finitely many points, and all we have to do is avoid these
intersection points in the path from $P_1$ to $Q_1$. We can also
assume that $\Delta t$ is small enough that $P_1$ crosses exactly
one boundary between times $t$ and $t+\Delta t$. Let us see how
$N_S$ changes in this time interval.

It will be convenient to call a circle $P_iP_jP_k$ {\it
$(a,b)$-splitting} (where $a+b = 2n-2$) if it has $a$ points of
$S$ inside it and the remaining $b$ outside it. For example, an
$(n-1,n-1)$-splitting circle is just a point-splitting circle.

\begin{figure}
\centering
\includegraphics[width=3in]{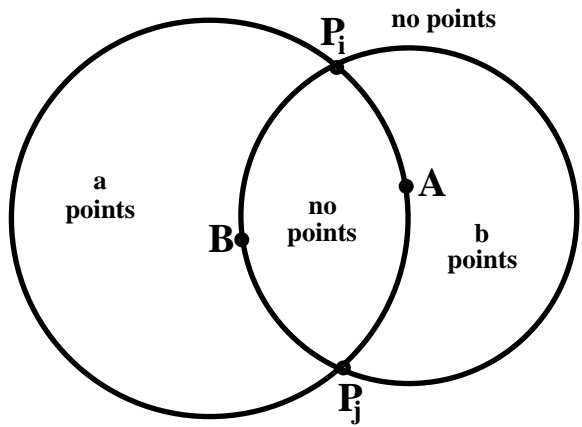}
\caption{$P_1$ crosses line $P_iP_j$.}
\end{figure}

First assume that $P_1$ crosses line $P_iP_j$, going from
position $P_1(t) = A$ to position $P_1(t+\Delta t) = B$. From the
remarks made above, we know that only circle $P_1P_iP_j$ can
change the value of $N_S$ by becoming or ceasing to be point-
splitting. Assume that circle $AP_iP_j$ was $(a,b)$-splitting.
Since $P_1$ only crossed the boundary $P_iP_j$ when going from $A$
to $B$, the region common to circles $AP_iP_j$ and $BP_iP_j$
cannot contain any points of $S$, as indicated in Figure 2. The
region outside of both circles cannot contain points of $S$
either. For circle $AP_iP_j$ to be $(a,b)$-splitting, the other
two regions must then contain $a$ and $b$ points respectively, as
shown. Therefore, circle $BP_iP_j$ is $(b,a)$-splitting. It
follows that $AP_iP_j$ was point-splitting if and only if
$BP_iP_j$ is point-splitting (if and only if $a=b=n-1$). We
conclude that the value of $N_S$ doesn't change when $P_1$
crosses a line determined by the other points; it can only change
when $P_1$ crosses a circle.

Now assume that $P_1$ crosses circle $P_iP_jP_k$, going from
position $P_1(t) = A$ inside the circle to position $P_1(t+\Delta
t) = B$ outside it. (The other case, when $P_1$ moves inside the
circle, is analogous.) We can assume that $P_1$ crossed the arc
$P_iP_j$ of the circle that doesn't contain point $P_k$. Notice
that $A$ must be outside triangle $P_iP_jP_k$ if we want $P_1$ to
cross only one boundary in the time interval considered. Assume
that circle $AP_iP_j$ was $(a,b)$-splitting. As before, we know
that the only regions of Figure 3 containing points of $S$ are
the one common to circles $AP_iP_j$ and $BP_iP_j$, and the one
outside both of them. They must contain $a-1$ and $b$ points
respectively, for circle $AP_iP_j$ to be $(a,b)$-splitting. In
this case, the value of $N_S$ can change only by circles
$P_iP_jP_k$, $P_1P_jP_k$, $P_1P_kP_i$ and $P_1P_iP_j$ becoming or
ceasing to be point-splitting. It is clear that circle
$P_iP_jP_k$ went from being $(a,b)$-splitting to being
$(a-1,b+1)$-splitting. The same is true of circle $P_1P_iP_j$.

\begin{figure}
\centering
\includegraphics[width=4in]{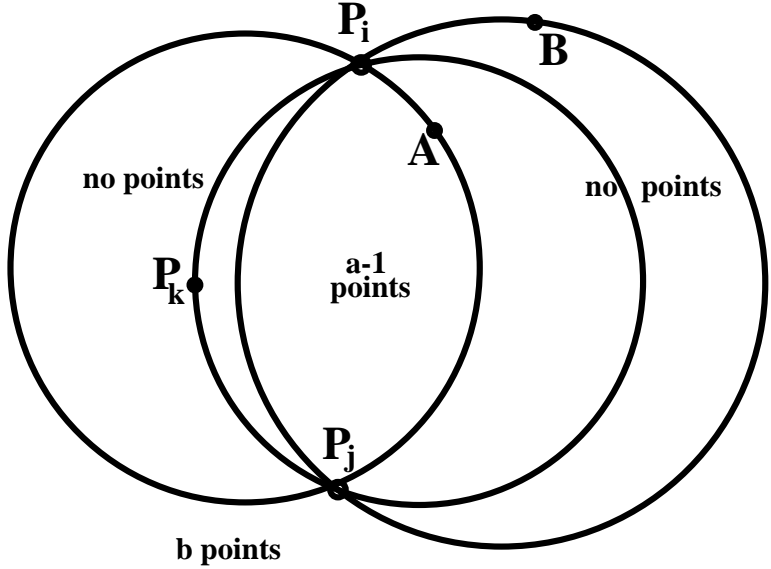}
\caption{$P_1$ crosses circle $P_iP_jP_k$.}
\end{figure}

It is also not hard to see, by a similar argument, that circles
$P_1P_jP_k$ and $P_1P_kP_i$ both went from being
$(a-1,b+1)$-splitting to being $(a,b)$-splitting. Again, the key
assumption is that $P_1$ only crossed the boundary $P_iP_jP_k$ in
this time interval.

So, by having $P_1$ cross circle $P_iP_jP_k$, we have traded two
$(a,b)$-splitting and two $(a-1,b+1)$-splitting circles for two
$(a-1,b+1)$-splitting and two $(a,b)$-splitting circles,
respectively. It follows that the number $N_S$ of point-splitting
circles remains constant when $P_1$ crosses a circle $P_iP_jP_k$
also.

We had shown that, as we moved $P_1$ to $Q_1$, $N_S$ could only
possibly change in a time interval when $P_1$ crossed a boundary
determined by the other points. But now we see that, even in such
a time interval, $N_S$ does not change! Therefore moving $P_1$ to
$Q_1$ doesn't change the value of $N_S$. Similarly, moving $P_i$
to $Q_i$ doesn't change $N_S$ either, for any $1 \leq i \leq
2n+1$. It follows that $N_S$ is the same for $S_{min}$ and
$S_{max}$. In fact, $N_S$ is the same for any set $S$ of $2n+1$
points in general position!

\medskip

\section {$N_S = n^2$}

\medskip

Now that we know that the number $N_S$ depends only on the number
of points in $S$, define $N_{2n+1}$ to be the number of
point-splitting circles for a set of $2n+1$ points in general
position. We compute $N_{2n+1}$ recursively.

\medskip

Construct a set $S$ of $2n+1$ points as follows. First consider
the vertices of a regular $2n-1$-gon with center $O$. Now move
them very slightly to positions $P_1, \ldots, P_{2n-1}$ so that
they are in general position. The difference will be so slight
that all the lines $OP_i$ still split the remaining points into
two sets of equal size, and all the circles $P_iP_jP_k$ still
contain $O$. Also consider a point $Q$ which is so far away from
the others that it lies outside of all the circles formed by the
points considered so far. Of course, we need $Q$ to be in general
position with respect to the remaining points. Let us count the
number of point-splitting circles of $S=\{O, P_1, \ldots,
P_{2n-1}, Q\}$.

First consider the circles of the form $P_iP_jP_k$. These circles
contain $O$ and don't contain $Q$; so they are point-splitting
for $S$ if and only if they are point-splitting for $\{P_1,
\ldots, P_{2n-1}\}$. Thus there are $N_{2n-1}$ such circles.

Next consider the circles $OP_iP_j$. It is clear that these
circles contain at most $n-2$ other $P_k$'s. They do not contain
$Q$, so they contain at most $n-2$ points, and they are not
point-splitting.

Finally consider the circles that go through $Q$ and two other
points $X$ and $Y$ of $S$. Circle $QXY$ splits the remaining
points in the same way that line $XY$ does. More specifically,
circle $QXY$ contains a point $P$ of $S$ if and only if $P$ is on
the same side of line $XY$ that $Q$ is. This follows easily from
the fact that $Q$ lies outside circle $PXY$. Therefore we have to
determine how many lines determined by two of the points of
$S-\{Q\}$ split the remaining points of this set into two sets of
$n-1$ points each. This question is much easier; it is clear from
our construction that the lines $OP_i$ do this and the lines
$P_iP_j$ do not. It follows that the $2n-1$ circles $OP_iQ$ are
point-splitting, and the circles $P_iP_jQ$ are not.

Summarizing, the point-splitting circles of $S$ are the
$N_{2n-1}$ point-splitting circles of $S-\{O,Q\}$ and the $2n-1$
circles $OP_iQ$. Therefore $N_{2n+1}=N_{2n-1}+2n-1$. Since
$N_3=1$, it follows inductively that $N_{2n+1}=n^2.$ This
completes the proof of Theorem 1.3.

\medskip

{\bf Theorem 3.1.} Consider a set of $2n+1$ points in general
position in the plane, and two  non-negative integers $a<b$ such
that $a+b=2n-2$. There are exactly $2(a+1)(b+1)$ circles which
are either $(a,b)$-splitting or $(b,a)$-splitting for the set of
points.

\medskip

{\it Sketch of Proof.} The argument of Section 2 carries directly
to this situation, to show that the number of circles in
consideration, which we denote $N(a,b)$, only depends on $a$ and
$b$. Therefore it suffices to compute it recursively, using the
set $S$ above. It is essential in the proof that $a<n-1$.

Just as above, there are $N(a-1,b-1)$ such circles among the
circles $P_iP_jP_k$. Among the $OP_iP_j$ there are exactly $2n-1$
such circles, namely the circles $OP_iP_{i+a+1}$ (taking
subscripts modulo $2n-1$). There are also $2n-1$ such circles
among the $QP_iP_j$, namely the circles $QP_iP_{i+a+1}$. Finally,
there are no such circles among the $OP_iQ$. Therefore $N(a,b) =
N(a-1,b-1) + 4n - 2 = N(a-1,b-1) + 2a + 2b + 2$.

Repeating the above argument for $a=0$, we get that
$N(0,b)=2b+2$. If we combine this and the recursive relation
obtained, Theorem 3.1 follows by induction.

\medskip

It is worth mentioning at this point that Theorems 1.3 and 3.1 are
closely related to a beautiful result of D.T. Lee, which gives a
sharp bound for the number of vertices of an order $j$ Voronoi
diagram. The connection is obtained if we embed our set $S$ of
points on the surface of a sphere. Then the point-splitting
circles of $S$ are put in correspondence with the
``point-splitting planes" of a three-dimensional convex polytope
with $2n+1$ vertices. These are known to be related to Voronoi
diagrams. See \cite[p. 397]{1} for more details on this, and a
proof of a result essentially equivalent to Theorems 1.3 and 3.1.

\medskip

\section{Questions}

Our work completely determines the number of point-splitting
circles, as well as the total number of $(a,b)$-splitting and
$(b,a)$-splitting circles for a set of points in general position
in the plane. However, we know very little about these numbers
for sets of points that are not in general position.

The situation here is much more subtle. For example, the number
of point-splitting circles of a set $S$ is not uniquely
determined by the subsets of $S$ which are concyclic. Consider
the following example. Let $S_1$ and $S_2$ be the two sets of
seven points shown in Figure 4. Both of them are almost in
general position; the only exception is that, for each of the two
sets, there is a circle going through four points of the set. In
$S_1$, this circle $\Gamma_1$ contains exactly one point of $S_1$
inside it. In $S_2$, this circle $\Gamma_2$ contains no points of
$S_2$ inside it. In analogy with Theorem 1.3, where it did not
matter which points were inside which circles, we might hope that
$S_1$ and $S_2$ have the same number of point-splitting circles.

\begin{figure}
\centering
\includegraphics[width=4in]{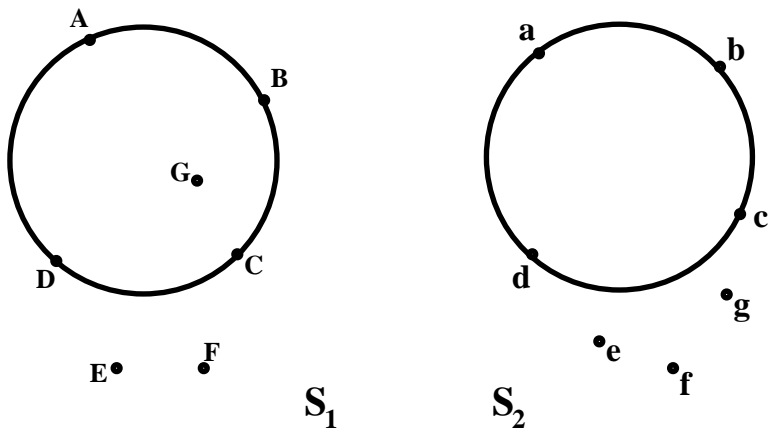}
\caption{$N_{S_1} = 8$ and $N_{S_2} = 9$.}
\end{figure}

Unfortunately this is not the case. If we move $A,B,C$ or $D$
very slightly to put $S_1$ in general position, the resulting set
will have $9$ point-splitting circles by Theorem 1.3. It is easy
to see that exactly two of $ABC, BCD, CDA$ and $DAB$ are among
these circles. When we deform the set back to $S_1$, these $9$
circles will still be point-splitting, but two of them will
deform into $\Gamma_1$. So $S_1$ has $8$ point-splitting circles.

Similarly, if we move $a,b,c$ or $d$ very slightly to put $S_2$
in general position, the resulting set will have $9$
point-splitting circles. But now we can see that when we deform
the set back to $S_2$, none of these circles will deform into
$\Gamma_2$, because $\Gamma_2$ contains no points of $S_2$.
Therefore $S_2$ has $9$ point-splitting circles.

Even if the number of point-splitting circles is not constant, we
might be able to say something about it. As a small example,
consider all sets of seven points which are almost in general
position, except that four of them are concyclic. It is possible
to show, by an argument similar to the above, that such a set can
only have $8$ or $9$ point-splitting circles. It seems reasonable
that, in general, one might be able to define some measure of how
far a set $S$ is from being in general position, and to obtain
bounds for $N_S$ in terms of that measure.

\medskip

\section{Acknowledgements}

The results in this article were first presented by the author at
Dan Kleitman's 65th birthday conference at the Massachusetts
Institute of Technology in 1999. The author wishes to thank Dan
Kleitman and Richard Stanley for useful discussions on this
subject during the conference.

He also wishes to thank Timothy Chow, who took an interest in
this result and mentioned it to several people. It is due to his
efforts that the author became aware of the equivalent result in
\cite{1} and the connection with Voronoi diagrams.

\end{document}